\documentclass[3p,times]{elsarticle}

\usepackage{ecrc,amssymb,amsmath,amsthm}

\volume{00}

\firstpage{1}

\runauth{}

\jid{procs}

\newcommand{\cP}{\mathcal{P}}
\newcommand{\mP}{\mathbb{P}}

\newcommand{\PGL}{\mathsf{PGL}}

\def\qed{{\hfill\hphantom{.}\nobreak\hfill$\Box$}}
\def\<{\langle}
\def\>{\rangle}

\newcommand{\cV}{\mathcal{V}}

\newcommand{\K}{\mathbb{K}}

\newtheorem{theorem}{Theorem}[section]
\newtheorem*{Theorem}{Main Result}
\newtheorem{remark}[theorem]{Remark}
\newtheorem{prop}[theorem]{Proposition}
\newtheorem{lemma}[theorem]{Lemma}
\newtheorem{observ}[theorem]{Observation}
\newtheorem{cor}[theorem]{Corollary}

\usepackage{amssymb}

\usepackage[figuresright]{rotating}

\begin{document}

\begin{frontmatter}



\dochead{}

\title{Severi's theorem for $d$-uple Veronese varieties\\Le th\'eor\`eme de Severi pour les vari\'et\'es $d$-uple de Veronese }

\author[JS]{Jeroen Schillewaert \fnref{JJ}}
\ead{jschillewaert@gmail.com}
\ead[JS]{www.jeroenschillewaert.com}
\fntext[JJ]{Research of the first author supported by Marie Curie IEF grant GELATI (EC grant nr 328178)}

\author[KS]{Koen Struyve \fnref{KK}}
\ead{kstruy@gmail.com}
\ead[KS]{http://cage.ugent.be/~kstruyve/}
\fntext[KK]{Research of the second author supported by the Fund for Scientific Research -- Flanders (FWO - Vlaanderen)}

\address[JS]{Department of Mathematics, Imperial College, London, United Kingdom}
\address[KS]{Department of Mathematics, Ghent University, Ghent, Belgium}


\begin{abstract}
We characterize $d$-uple Veronese embeddings of finite-dimensional projective spaces. The easiest non-trivial instance of our theorem is the embedding of the projective plane in 5-dimensional projective space, a result obtained in 1901 by Severi when the underlying field is complex.\vspace{1mm}    

\flushleft{\normalsize \bf Resum\'e}\\
Le plongement de Veronese $d$-uple d'un espace projectif de dimension finie est caracteris\'e. L'exemple le plus simple et non-trivial de notre th\'eor\`eme est le plongement 
d'un plan projectif dans un espace projectif de dimension 5, un r\'esultat obtenu par Severi en 1901 si le corps sous-jacent est le corps des nombres complexes.
\end{abstract}


\end{frontmatter}



\section{Introduction}
In 1901 Severi \cite{Severi} characterized the Veronese surface as the only smooth irreducible and non-degenerate projective surface in a 5-dimensional complex projective space which can be projected isomorphically into a 4-dimensional complex projective space.

In 1984, Mazzocca and Melone \cite{Maz-Mel:84} used three simple axioms (the MM-axioms) to characterize the ordinary quadric Veronese variety over finite fields of odd order. This was extended to all finite fields by Hirschfeld and Thas \cite{Hir-Tha:91} and to arbitrary fields by Schillewaert and Van Maldeghem \cite{SVM2}. Our main result extends Severi's original characterization \cite{Severi} to arbitrary finite degree $d$ and arbitrary fields of cardinality at least $d^{\frac{3}{2}}$ using a generalization of the MM-axioms, as Severi's conditions imply the MM-axioms \cite{LVDV,Zak}.

\subsection{Definitions}\label{section:defs}

The \emph{Veronese embedding} $\phi_{n,d}:\mathbb{P}^n(\K)\to \mathbb{P}^{{n+d \choose d}-1}(\K)$ maps $(x_0,\cdots,x_n)$ to $(x_0^d,x_0^{d-1}x_1,\cdots,x_n^d)$. The \emph{Veronese variety} $\mathcal{V}_n^{d^n}$ is the image of this map. A \emph{rational normal curve} $C$ in a $d$-dimensional projective space $\Sigma$ is a set of points of $\Sigma$ which is projectively equivalent to $\mathcal{V}_1^{d}$.  For any $x \in C$ one can define a unique {\em tangent line} $T_x(C)$ which is determined by the points of the curve if either $d=2$ or $d \geq | \K |+2$.


Let $V$ be a possibly infinite-dimensional, non-trivial vector space over some field $\K$, and let $\mathbb{P}(V)$ be the corresponding projective space.  Let $X$ be a spanning point set of $\mathbb{P}(V)$ and let $\Xi$ be a collection of  $d$-dimensional subspaces of $\mathbb{P}(V)$, the {\em rational spaces} of $X$, such that, for any $\xi\in\Xi$, the intersection $\xi\cap X$ is a rational normal curve $X(\xi)$ in $\xi$ of dimension $d$  ($d \geq 2$) and then, for $x\in X(\xi)$, we sometimes denote $T_x(X(\xi))$ simply by $T_x(\xi)$. We call $\mathcal{V}=(X,\Xi)$, or briefly $X$, a {\em Veronesean cap of degree $d$} if the following properties (V1), (V2) and (V3) hold.
\begin{itemize}\addtolength{\itemsep}{-0.4\baselineskip}

\item[(V1)] Any two points $x$ and $y$ of $X$ lie in a unique element of $\Xi$, denoted by $[x,y]$.

\item[(V2)] If $\xi_1,\xi_2\in \Xi$, with $\xi_1\neq \xi_2$, then $\xi_1\cap\xi_2\subset X$.

\item[(V3)] For every $x\in X$, every $\xi\in\Xi$, with $x\notin \xi$, $\cup_{y\in X(\xi)} T_x([x,y])$ is a plane $T(x,\xi)$  depending on $x$ and $\xi$.
\end{itemize}

%


In the complex case, every (possibly singular) non-degenerate curve in $\mathbb{P}^{d}$ of degree $d$ is a rational normal curve \cite[p.121]{Shafarevich}. In the finite case one can sometimes relax the conditions \cite{Casse-Glynn}, but not always \cite{Glynn,Segre}.

\subsection{A Veronese variety is a Veronesean cap}

By Lemma 2.3 of \cite{KantorShult} the images of lines are rational normal curves. Property (V1) is immediate. To verify (V2), by homogeneity, see \cite{KantorShult}, we only need to check two situations, namely $\psi_{1,2}:=\<\phi_{n,d}(L_1)\>\cap \<\phi_{n,d}(L_2)\>$ and if $n\geq 3$ we need to compute also $\psi_{1,3}:=\<\phi_{n,d}(L_1)\>\cap\<\phi_{n,d}(L_3)\>$ where $L_1$ is given by $X_i=0,2\leq i\leq N$, $L_2$ by $X_i=0, i\neq 2,1\leq i\leq N $ and $L_3$ by $X_i=0;0\leq i \leq N-2$. Then a computation yields $\psi_{1,2}:=\phi_{n,d}(L_1\cap L_2)$ and $\psi_{1,3}=\emptyset$. 
For (V3), again by homogeneity, consider the point $r:=(0,\dots,0,1)$ and the points $p_1:=(1,0,\dots,0)$ and $p_2:=(0,1,0,\dots,0)$ on $L_1$.
Then a computation yields that the plane spanned by the tangent lines at $\phi_{n,d}(r)$ to $\phi_{n,d}(\langle r,p_1 \rangle)$ and $\phi_{n,d}(\langle r,p_2 \rangle)$
is the requested one.


\section{Statement of the main results}

\begin{prop}\label{tangentspace}
If $X$ is a Veronesean cap, then for every $x\in X$ the set $\cup_{\xi \in\Xi, x \in \xi } T_x(\xi)$ is a subspace of constant dimension.
\end{prop}

We call the dimension of these subspaces the \emph{index}. The next result classifies Veronesean caps of sufficiently large dimension with respect to its degree and index over fields of sufficiently large cardinality with respect to the degree.

\begin{Theorem}\label{Thm:main}
If $X$ is a Veronesean cap of degree $d$ and finite index $n\geq 2$ in a projective space $\mP(V)$ over a field $\K$ such that $\dim V\geq M:= {n+d \choose d}$ and such that $|\K | \geq d^{3/2}$, then $\dim V=M$, and $X$ is projectively equivalent to the Veronese variety $\mathcal{V}_n^{d^n}$ over $\K$.
\end{Theorem}

\begin{remark}\label{rem:ineq}
One can strengthen the bound on the field size in the Main Result to the combination of the bounds $$|\K | \geq d + 2   \mbox{ when $d \neq 2$, and } \frac{|\K|^{n} - 1}{|\K| -1} \geq 2 + \sum_{i=2}^{\min(n-1,d) } (-1)^{i} {d+1 \choose i+1} \frac{|\K|^{n-i} -1}{|\K|-1}.$$ The first bound is necessary as otherwise the tangent lines are not determined by the points of the curve. The second bound becomes irrelevant in view of the first one for $n \leq 2$ or $d \leq 10$.
\end{remark}


\section{Basic properties of rational normal curves}\label{rncprop}
Let $C$ be a rational normal curve in a $d$-dimensional projective space $\Sigma$ as defined in Section~\ref{section:defs} (with $d \geq 2$).

\begin{lemma}
Every $d+1$ points are linearly independent. \cite[Theorem 1.1]{KantorShult}
\end{lemma}

\begin{observ}\label{obs:1}
The rational normal curve $C$ comes equipped with a notion of crossratio on quadruples of points, such that the group $\PGL(2,\K)$ acts three-transitively on $C$ while preserving the crossratio. (See~\cite[Ex. 1.19]{Har:92})
\end{observ}

\begin{lemma}\label{lem:proj1}
Let $x$ be a point of the curve $C$, then a projection from $x$ to a hyperplane not incident with it maps $C \setminus \{x\}$ and the tangent line through $x$ bijectively to the points of a rational normal curve in a $(d-1)$-dimensional projective space. This projection preserves the crossratio on $C$.
\end{lemma}
\proof
By Observation~\ref{obs:1} we may fix $x=(1,0,\dots,0)$ in $\Sigma$. This reduces the problem to an easy calculation.\qed

Repeated application of Lemma~\ref{lem:proj1} shows
\begin{cor}\label{cor:proj1}
Let $x_1, \dots, x_i$ be $i$ pairwise different points of $C$ such that $1 \leq i \leq d-1$. Then a projection from the span of $x_1, \dots, x_i$ to a complementary subspace maps $C \setminus \{x\}$ and the tangent lines through $x_1, \dots, x_i$ bijectively to the points of a rational normal curve in a $(d-i)$-dimensional projective space. This projection preserves the crossratios $C$.
\end{cor}

\begin{observ}\label{obs:2}
By considering 4 points of $C$ and applying Corollary~\ref{cor:proj1} to $d-1$  other points on the curve (which always exist when $| \K | \geq 4$ by the first bound of Remark~\ref{rem:ineq}), one obtains that the crossratio mentioned in Observation~\ref{obs:1} solely depends on the point set of $C$. For $|\K| = 2$ or $3$ this is trivial as there is only one way to define a valid crossratio on 3 or 4 points.
\end{observ}

Using Corollary~\ref{cor:proj1} considering a projection from $\langle C\cap L_1,\cdots,C\cap L_i \rangle$ yields
\begin{cor}\label{cor:proj2}
Let $L_1, \dots, L_i$ and $x_1, \dots x_j$ be respectively $i$ tangent lines on $C$ and $j$ points of $C$, none of which are equal or incident with each other. If $2i +j \leq d+1$, then the span of these lines and points is of dimension $2i +j -1$.
\end{cor}

\section{Proof of the Main Result}

\subsection{The projective space associated with a Veronesean cap}\label{projspace}
We use the same notation as in Section~\ref{section:defs}. Associated with $\cV$ we can consider the geometry $\cP$ having point set $X$ and line set $\Xi$, endowed with the natural incidence.
Then a proof similar to~\cite[Prop 2.2]{ThasHVM} yields
\begin{prop}\label{projplane}
$\cP$ is  a projective space of dimension $n \geq 2$.
\end{prop}

%
%

Also the proof of Proposition~\ref{tangentspace} is essentially the same as the proof of~\cite[Prop. 2.4]{ThasHVM}. The index obtained is the same $n$ as in Proposition~\ref{projplane}.

From now on, we denote a point of $\cP$ with a $\bar{.}$, e.g. $\bar{x}$, and the corresponding point in $\mP(V)$ without a $\bar{.}$, e.g. $x$. Similarly we denote subspaces of $\cP$ with a $\bar{.}$, and if $\bar\pi$ is a subspace of $\cP$, then $\pi:=\langle  x \in X \vert \bar{x} \in \bar\pi  \rangle$. The following easy observation will be crucial for induction arguments.

\begin{lemma}\label{lem:ind}
Let $\bar\pi$ be some non-empty subspace of $\cP$, then $(X_\pi, \Xi_\pi)$ with $X_\pi := \{x \vert \bar{x} \in \bar{\pi}\}$, and $\Xi_\pi :=\{\xi \in \Xi \vert  X(\xi) \subset X_\pi\}$, is a Veronesean cap of degree $d$ inside $\pi$, which we call a {\em subcap}.
\end{lemma}

\subsection{Dimensional analysis}

In the remainder of the proof let $(X,\Xi)$ be a Veronesean cap of degree $d$ and index $n$, satisfying the restrictions listed in the Main Result and let $M={n+d \choose d}$. 
We assume as induction hypothesis that the Main Result holds for Veronesean caps of index up to $n-1$,  note that the cases of index 0 and 1 are trivial.

\begin{lemma}\label{lem:secants}
Let $\bar\pi_1, \bar\pi_2, \dots \bar\pi_{d+1}$ be a set of hyperplanes of $\cP$ in general position. If $\bar{x}$ is a point of $\cP$ not contained in any of these hyperplanes, then there exist at least two different lines through $x$ such that each point on these lines is contained in at most two of the hyperplanes  $\bar\pi_1, \bar\pi_2, \dots \bar\pi_{d+1}$.
\end{lemma}
\proof
This holds as the right hand side of the second inequality of Remark~\ref{rem:ineq} counts the number of intersection points of at least three of these subspaces (via the inclusion-exclusion principle) plus two, while the left hand side of this inequality counts the number of lines through $x$.
\qed

For the following proposition, note that whenever $| \K | \geq d$ there do exist hyperplanes $\bar\pi_1, \bar\pi_2, \dots, \bar\pi_{d+1}$ in general position, for example the dual of the point set of a rational normal curve in $\cP$.

\begin{prop}\label{prop:bound}
Let $\bar\pi_1, \bar\pi_2, \dots, \bar\pi_i$ ($0 \leq i \leq d+1$) be a number of hyperplanes of $\cP$ in general position which is extendable to a set of $d+1$ hyperplanes in general position. Then
\begin{itemize}\addtolength{\itemsep}{-0.4\baselineskip}
\item[(i)] $\mathrm{codim} \<\pi_1, \pi_2 , \dots, \pi_i \> ={n + d  -i\choose d -i}$.
\item[(ii)]$\dim V=M$. 
\item[(iii)] $\pi_i \cap \<\pi_1, \dots,\pi_{i-1}  \>=\< \pi_i \cap \pi_1, \dots, \pi_i \cap \pi_{i-1}  \>$.
\item[(iv)]  If the geometric dimension of $\bar\pi$ is $m$, then the geometric dimension of  $\pi$ is ${m+d \choose d} -1$. In particular  $(X_\pi, \Xi_\pi)$ is a Veronesean cap of degree $d$ and index $m$ satisfying the restrictions listed in the Main Result.
\end{itemize}
\end{prop}
\proof 
First we prove $\mathrm{codim} \<\pi_1, \pi_2 , \dots, \pi_i \>  \leq {n + d  -i\choose d -i}$, not using $\dim V \geq M$, by induction on both $n$ and $d-i$, making use of Lemma~\ref{lem:ind}.  The case $n = 1$ follows by Section \ref{rncprop}, so assume $n\geq 2$.

First assume that $d- i = -1$. Let $x$ be any point in $X \setminus (\cup_{j=1}^{d+1} \pi_j  )$, and let $\xi \in \Xi$ be a rational normal space through $x$ such that each $y \in X(\xi)$ is contained in at most two of the subspaces $\pi_i$ (which exists by Lemma~\ref{lem:secants}).
For such a $\xi$, we have that if $y$ is contained in some subspace $\pi_i$, then it is contained in $\<\pi_1, \pi_2 , \dots \pi_i \>$. If it is contained in two of these subspaces, then also its tangent line is contained in $\<\pi_1, \pi_2 , \dots \pi_{d+1} \>$. Hence $\<\pi_1, \pi_2 , \dots \pi_{d+1} \>$ contains $\xi$, and by varying $x$ all of $X$.

Secondly assume that $d-i \geq 0$. In this case we pick a hyperplane $\bar\pi_{i+1}$ of $\cP$ such that $\bar\pi_1, \bar\pi_2, \dots, \bar\pi_{i+1}$ are in general position. The codimension of the span $ \<\pi_1, \pi_2 , \dots ,\pi_i \> $ in $\mP(V)$ equals the sum of the codimension $c_1$ of $ \<\pi_1, \pi_2 , \dots, \pi_i \> $ in $ \<\pi_1, \pi_2 , \dots, \pi_{i+1} \> $ and the codimension $c_2$ of $ \<\pi_1, \pi_2 , \dots, \pi_{i+1} \> $ in $\mP(V)$. Then $c_1$ is at most the codimension of $ \<\pi_1\cap \pi_{i+1} , \dots , \pi_i \cap \pi_{i+1}\> $ in $\pi_{i+1}$, which is bounded from above by ${(n-1)+d - i  \choose  d-i }$ by induction on $n$, and $c_2$ is bounded from above by ${n+d - (i + 1) \choose  d-(i+1) }$ by induction on $d-i$. The sum of both is at most ${ n+d - i \choose d-i }$.

Combining the bound for $i=0$ with the restriction $\dim V \geq M$ yields $\dim V = M$ (proving (ii)), and that all bounds above have to be sharp, which proves (i) and (iii). Repeated application of (i) with $i=1$ proves (iv).\qed

%
%
%

\subsection{The projective space $\cP$ is isomorphic to $\mathbb{P}^{n}(\K)$} \label{section:projspace}

As we know the dimension of $\cP$  from Proposition~\ref{projplane} one can use part (iv) of Proposition~\ref{prop:bound} to restrict to the case $n=2$ for the remainder of Section~\ref{section:projspace}. 
Remember from Observations~\ref{obs:1} and~\ref{obs:2} that the point sets of the rational normal curves $\xi \cap X$ with $\xi \in \Xi$ come equipped with a natural crossratio on quadruples. By construction of $\cP$ this yields a crossratio on the point sets of lines of $\cP$.  The next step is now to define a crossratio on line pencils of $\cP$.

\begin{observ}
Let  $\bar{x}$ be an arbitrary point of $\cP$ and identify each line $\bar\zeta$ through $\bar{x}$ with the tangent line $T_x(\zeta)$. Hence, by (V3), we obtain a correspondence of lines through $\bar{x}$ with lines in the tangent plane on $x$. The latter, being a line pencil in the projective space $\mP(V)$, comes with a natural notion of a crossratio, which via identification provides a crossratio on the line pencil through $\bar{x}$ in $\cP$.
\end{observ}

The next lemma links these crossratios together.

\begin{lemma}\label{lem:cros} The perspectivity from the line pencil through a point $\bar{x}$ of $\cP$ to the points on a line $\bar\xi$ of $\cP$ not containing $\bar{x}$ preserves the crossratios defined on these.
\end{lemma}
\proof
Extend the line $\bar\xi$ to a set of lines $\{\bar\xi, \bar\xi_2 , \dots \bar\xi_{d+1}  \}$ in general position, such that only $\bar\xi_{d+1}$ contains $\bar{x}$. 

Parts (i) and (ii) of Proposition~\ref{prop:bound} imply that the codimension of $\langle  \xi_2, \dots, \xi_{d} \rangle$ is 3. Hence we may consider a projection $\tau$ of $\mP(V)$ from $\langle \xi_2, \dots, \xi_{d} \rangle$ on some two-dimensional subspace disjoint from it. Proposition~\ref{prop:bound} also implies that $\tau$ maps $\xi$ on some line $L$ (as $\langle \xi, \xi_2, \dots, \xi_{d} \rangle$ is of codimension one)  and that  $x$ is mapped to a point disjoint from $L$ (by part (iii) applied to $\xi_{d+1}$ and $\langle \xi, \xi_2, \dots, \xi_{d} \rangle$).  

The points of the curve $X(\xi)$ not in $ \xi_2 \cup \dots \cup \xi_{d}$ and the tangent lines on $\xi$ through the remaining points on this curve project to pairwise distinct points of $L$ while preserving the crossratios,  by Corollary~\ref{cor:proj1} and part (iii) of Proposition~\ref{prop:bound}. This Corollary~\ref{cor:proj1} also states that this projection preserves crossratios, hence we have connected the crossratio on $\xi \cap X$ and the crossratio on $L$.

Let $\xi'$ be some rational normal space through $x$. From Corollary~\ref{cor:proj2} and counting the number of  points of $X \cap \xi'$ contained in any of the $\xi, \xi_2 , \dots \xi_{d-1}$, while accounting for the tangent lines of the points that are contained in two of these rational normal spaces (applying (V3)), it follows that the image of $\xi'$ under $\tau$ is at most one-dimensional. As this image contains $\tau(x)$ and a point of $L$ it has to be a line.

In the image of the projection we hence obtain a perspectivity between the line pencil through the point $\bar{x}$ of $\cP$ and the points on $\bar\xi$, preserving the crossratios on both. \qed

This has the following immediate consequence.

\begin{cor}
The induced action on the rational normal curve $X(\xi)$ of the projectivity group of $\cP$ w.r.t. to a line $\bar\xi$ of $\cP$ preserves the natural crossratio on it.
\end{cor}

Let $G$ be the projectivity group of $\cP$ w.r.t. to some line $\bar\xi$ of $\cP$.  Then $G \leq \mathrm{PGL}(2,\K)$ as the latter consists of those permutations preserving the crossratio. As $\mathrm{PGL}(2,\K)$ acts sharply 3-transitive while $G$ is at least 3-transitive, one has $G=\mathrm{PGL}(2,\K)$. Moreover, by von Staudt's theorem \cite{VonStaudt}, this implies that the projective plane $\cP$ is Pappian, hence $\K$ is a field and $\cP$ is isomorphic to $\mathbb{P}^{2}(\K)$.

Let us mention another corollary.

\begin{cor}\label{cor:findtangents}
Let $\bar{x}$ and $\bar\xi$ be a point and line of $\cP$ such that $\bar{x} \notin \bar\xi$. Then the map 
\begin{align*}
X(\xi)  \to T(x,\xi) : y  \mapsto  T_x([x,y])
\end{align*} 
is completely determined by the curve $X(\xi)$, the point $x$, and the images of at least three points in $X(\xi)$ under this map.
\end{cor}
\proof
By Lemma~\ref{lem:cros} this map preserves crossratios, so  knowing the image of three points suffices to completely determine the map.
\qed

\subsection{Mapping $X$ onto a Veronese variety}

Let $Y$ represent the point set of the projective space $\cP$. By construction of $\cP$ we have a bijective map $\iota: Y \to X: \bar{x} \mapsto x$. If we set $X'$ and $\Xi'$ to be respectively the point set and the set of rational spaces of a Veronese variety $\mathcal{V}^{d^n}_n$ over the field $\K$ with ambient vector space $\mP(V')$, then, by construction of the variety, we also have a bijective map $\iota' : Y \to X'$, i.e. the Veronese embedding.

Our goal is to show that the bijective map $\varphi := \iota' \circ \iota^{-1}$ from $X$ to $X'$ lifts uniquely to a collineation between the ambient projective spaces. We do this by induction on subcaps $Z$ of increasing index $m$. 

%

\subsubsection{The case $m=1$}
Here we have to consider a map $\varphi$ defined on a subcap $Z$ of index 1, so on the points of a certain rational normal curve $X(\xi)$ (with $\xi \in \Xi$). We want to show that $\varphi$ extends uniquely to a collineation defined on the subspace $\xi$ of $\mP(V)$. This is the case if and only if $\varphi$ preserves the crossratios on the rational normal curve up to a possible field automorphism (see for example~\cite[Ex. 1.19]{Har:92}).

We may assume that the index of $X$ is 2 (by using (iv) of Proposition~\ref{prop:bound}), so  that $\cP$ is a projective plane.
As the projectivity group determines crossratios up to a possible field automorphism $\sigma$, and by the fact that this group is a feature of $\cP$ and does not depend on the actual structure of $X$ and $X'$, it is implied that $\varphi$ preserves the crossratio on $X(\xi)$ up to a field automorphism $\sigma$.

\begin{remark} \label{rem:indep}
Moreover, by Lemma~\ref{lem:cros}, the field automorphism $\sigma$ does not depend on the choice of the subcap $Z$.
\end{remark}

\subsubsection{The case $m \geq 2$}
We may suppose without loss of generality that $Z=X$, so $m=n$. Fix $d+1$ hyperplanes $\bar\pi_1, \bar\pi_2, \dots, \bar\pi_{d+1}$ of $\cP$ in general position, obtained by considering the dual of a rational normal curve in $\mP^n(\K)$. The spaces  $\pi_1, \dots, \pi_{d+1}$ span $\mP(V)$ by Proposition~\ref{prop:bound}.

The induction hypothesis allows us to define a collineation $\rho$ mapping the subspace $\pi_1$ to the subspace $\varphi(\pi_1)$, in such a way that $\rho$ extends $\varphi$ on $X \cap \pi_1$. We now extend $\rho$ recursively to a collineation agreeing with $\varphi$ on $X \cap (\cup_{j=1}^{d+1} \pi_j)$. Suppose that we already have defined $\rho$ on the span of $\pi_1$ up to $\pi_i$, agreeing with $\varphi$ on $X \cap (\cup_{j=1}^{i} \pi_j)$ ($1 \leq i \leq d$). Then, in order to extend $\rho$ to a collineation defined on the span up to $\pi_1$ up to $\pi_{i+1}$, it suffices to extend the restriction of $\rho$ on $\pi_{i+1} \cap \<\pi_1, \dots , \pi_i \>$ to $\pi_{i+1}$. This is possible because of the induction hypothesis and (iii) and (iv) of Proposition~\ref{prop:bound} (and using Remark~\ref{rem:indep} in the case $m=2$ to avoid having conflicting field automorphisms). Eventually we obtain a collineation $\rho$ defined on the entirety of $\mP(V)$, agreeing with $\varphi$ on $X \cap (\cup_{j=1}^{d+1} \pi_j)$.

\begin{prop}
The collineation $\rho$ agrees with $\varphi$ on the entirety of  $X$. 
\end{prop}
\proof
Let $x$ be any point in $X \setminus \cup_{j=1}^{d+1} \pi_j$. It suffices to show that $x$ is uniquely determined by the points  $X \cap (\cup_{j=1}^{d+1} \pi_j)$. 

By Lemma~\ref{lem:secants} we can find two lines $\bar\xi_1$ and $\bar\xi_2$ through $\bar{x}$ such that each point on these lines is contained in at most two of the hyperplanes $\bar\pi_1, \bar\pi_2, \dots, \bar\pi_{d+1}$. If $d\leq 3$ we can pick these such that each point on these lines is contained in at most one such hyperplane. 

Assume $d > 3$, and let $y$ be a point on $\xi_1$ contained in two of the hyperplanes, we may suppose these are $\pi_1$ and $\pi_{d+1}$. By the bounds in Remark~\ref{rem:ineq} we have $| \K | \geq d + 2$. Next we reconstruct the tangent line $L:= T_y(\xi_1)$ from the points in  $X \cap (\cup_{j=1}^{d+1} \pi_j)$. We already know the position of $y$ and some plane in which $L$ is contained by property (V3). If $n >2$ we know at least two of these planes which determines $L$ uniquely, so we are left with the case that $n=2$. 

As  $| \K | \geq d+ 2$ and by the construction of the lines $\bar\pi_1, \bar\pi_2, \dots, \bar\pi_{d+1}$, there exists some lines $\bar\pi_{d+2}$ and $\bar\pi_{d+3}$ such that $\bar\pi_1, \bar\pi_2, \dots, \bar\pi_{d+3}$ are in general position. As $\bar\pi_{d+2}$ and $\bar\pi_{d+3}$ intersect the $\bar\pi_1, \bar\pi_2, \dots, \bar\pi_{d+1}$ in $d+1$ pairwise different points, the position of the subspaces $\pi_{d+2}$ and $\pi_{d+3}$ is uniquely determined. Also the unique intersection point of both, which we denote by $z$, is  determined. 

Applying part (iii) of Proposition~\ref{prop:bound} to  $\pi_1$ and $\pi_{d+2}$ yields that $y,z\notin \langle \pi_{i}\vert 2 \leq i\leq d \rangle$. Part (iii) of Proposition~\ref{prop:bound} applied to $\pi_1$ and $\pi_2, \dots, \pi_{d+1}$ also yields $T_y(\pi_1)\not\subset \langle \pi_{i} \vert 2 \leq i\leq d \rangle$, in particular the projections of the tangent lines $T_y(\pi_1)$ and $T_y(\pi_{d+1})$, and hence all tangent lines in $T(x,\pi_2)$, are pairwise different. Using this projection we can also determine the tangent line $T_y([y,z])$ (this is the unique line in $T(x,\pi_2)$ mapped to the line containing the projections of $y$ and $z$). Corollary~\ref{cor:findtangents} now allows us to determine all tangent lines through $y$, in particular $L$.

By Corollary~\ref{cor:proj2} the points in $X \cap (\cup_{j=1}^{d+1} \pi_j)$ uniquely determine $\xi_1$ and analogously $\xi_2$. The point $x \in \xi_1\cap\xi_2$ is hence also uniquely determined.
\qed

This concludes the proof of the Main Result.

\end{document}